\documentclass[11pt]{article}
\usepackage{amssymb}
\usepackage{color, amsmath,amssymb, amsfonts, amstext,amsthm, latexsym}

\usepackage{amssymb, epsfig, amssymb, latexsym}
\usepackage{amsmath}
\usepackage{graphicx}
\usepackage{longtable}

\newcommand{\s}{\sigma}

\renewcommand{\phi}{\varphi}

\renewcommand{\a}{\alpha}

\newcommand{\R}{{\mathbb R}}

\newcommand{\eps}{\varepsilon}

\newcommand{\PX}{{\Bbb{P}}}

\newtheorem{theorem}{Theorem}
\newtheorem{lemma}{Lemma}
\newtheorem{remark}{Remark}

\title{An intermediate regime for exit phenomena driven by
non-Gaussian L\'evy noises\footnote{This work was partly supported
by the NSF Grant 0620539. } }

\author{Zhihui Yang$^1$ and Jinqiao Duan$^2$\\
\\
1. Department of Mathematics\\Western Illinois University\\
Macomb, IL  61455, USA.\\
 \emph{E-mail: z-yang2@wiu.edu }
\\
2. Department of Applied Mathematics\\ Illinois Institute of Technology \\
  Chicago, IL 60616, USA.\\
 \emph{E-mail: duan@iit.edu}  }

\begin{document}
\date{May 18, 2008 (Revised)}

\maketitle

\pagestyle{plain}

\begin{abstract}
A dynamical system driven by non-Gaussian L\'evy noises of small
intensity   is considered. The first exit time of solution orbits
from a bounded neighborhood of an attracting equilibrium state is
estimated. For a class of non-Gaussian L\'evy noises, it is shown
that the mean exit time is asymptotically faster than exponential
(the well-known Gaussian Brownian noise case) but slower than
polynomial (the stable L\'evy noise case), in terms of the
reciprocal of the small noise intensity.

\medskip

%{\bf Short Title:} Exit time  under non-Gaussian L\'evy noises   \\



\bigskip

\emph{Dedicated to Professor Ludwig Arnold on the occasion of his
70th birthday}

\end{abstract}

\section{Introduction}  \label{intro}
%%%%%%%%%%%%%%%%%%%%%%%%%%%%%%%%%%%%%%%%%%%%%%%%%%%%%%%%%%%%%%%%%

Although Gaussian processes like Brownian motion have been widely
used in modeling fluctuations in   engineering and science, it turns
out that some complex phenomena  involve with non-Gaussian L\'evy
motions. For instance, it has been argued that diffusion by
geophysical turbulence \cite{Shlesinger} corresponds, loosely
speaking, to a series of  ``pauses", when the particle is trapped by
a coherent structure, and ``flights" or ``jumps" or other extreme
events, when the particle moves in the jet flow. Paleoclimatic data
\cite{Dit} also indicates such irregular processes.

L\'evy motions are thought to be appropriate models for non-Gaussian
processes with jumps \cite{Sato-99}. Let us recall that  a L\'evy
motion $L(t)$, or $L_t$, has  independent and stationary increments,
i.e., increments $\Delta L (t, \Delta t)= L(t + \Delta t)-  L(t)$
are stationary (therefore $\Delta L$ has no statistical dependence
on $t$) and independent for any non overlapping time lags $\Delta
t$. Moreover, its sample paths are only continuous in probability,
namely, $\PX (|L(t)-L(t_0)| \geq \delta) \to 0$ as $t\to t_0$ for
any positive $\delta$. This continuity is weaker than the usual
continuity in time.

This generalizes the Brownian motion $B(t)$, as $B(t)$ satisfies all
these three conditions.  But \emph{Additionally},  (i) Almost every
sample path of the Brownian motion     is  continuous in time in the
usual sense and (ii) Brownian motion's increments are Gaussian
distributed.

  SDEs  driven by   non-Gaussian L\'evy noises
have attracted much attention  recently \cite{Apple,  Schertzer}.
Although the SDEs driven by L\'evy motion may generate stochastic
flows \cite{Kunita2004, Apple}, or generate random dynamical
systems in the sense of Arnold \cite{Arnold}, under certain
conditions,
 this research issue is still under
development. Recently, mean exit time estimates have been
investigated by Imkeller and Pavlyukevich \cite{ImkellerP-06,
ImkellerP-08}.

\medskip

Consider a scalar deterministic ordinary differential equation
$\dot{Y}_{t}=-U^{'}(Y_{t})$, $Y_{0}=x \in [-b,a]$, $a,b >0$, where
the potential function $U$ is a sufficiently smooth   function.
Assume that $0$ is a asymptotically stable equilibrium. Namely, for
any starting point $x$ in $[-b,a]$, the trajectory   $Y_{t}$ tends
to $0$  as time $t \rightarrow \infty$. In this case, $U(\cdot)$ has
a minimum at $0$.

%with a minimum at the origin. Assume $U^{'}(x)
%>0$, $U(0)=0$, $U^{'}(x)=0$ if and only if  $x=0$, and $U^{''}(0)
%=M>0$. These assumptions imply that

\vskip 12pt

Now perturb the deterministic dynamical system $Y_{t}$ with some
small random noise. Let us consider a scalar stochastic differential
equation (SDE)
\begin{equation} \label{sde}
dX_{t}^{\varepsilon} =  - U^{'} (X_{t}^{\varepsilon}) dt
+\varepsilon d L_{t}, \;\; X_0 = x,
\end{equation}
 where  $0<\varepsilon \ll 1 $ is the noise intensity, and $L_{t}$ is
a L\'evy process. A scalar L\'evy process is characterized by a
drift parameter $\theta$, a variance parameter $d>0$ and   a non-negative
Borel measure $\nu$, defined on $(\R, \mathcal{B}(\R))$ and
concentrated on $\R \setminus\{0\}$, which satisfies
\begin{equation} \label{levycondition}
  \int_{\R \setminus\{0\} } (y^2 \wedge 1) \; \nu(dy) < \infty,
\end{equation}
or equivalently
\begin{equation}
  \int_{\R \setminus\{0\} } \frac{y^2}{1+y^2}\; \nu(dy) < \infty.
\end{equation}
This measure $\nu$ is the so called L\'evy measure or the L\'evy
jump measure of the L\'evy process $L(t)$. We also call $(\theta, d,
\nu)$ the \emph{generating triplet}.

\vskip 12pt

 We   study the first
exit problem for the solution process $X_{t}^{\varepsilon}$ from
bounded intervals containing the attracting equilibrium $0$, as
$\varepsilon \downarrow 0$.

We define the first exit time from  the spatial interval $[-b,a]$,
$a$ and $b$ positive, as follows:
\[
\sigma(\varepsilon) = \inf \{t \geq 0, X_{t}^{\varepsilon} \notin
[-b,a] \}.
\]

 It is known \cite{FW}
that when the noise is Gaussian Brownian, i.e., the L\'evy measure
part is absent,  the mean exit time of the perturbed system
$X_{t}^{\varepsilon}$ is asymptotically exponentially fast
\begin{equation}
E_x \s(\eps) \sim \exp(\frac{C}{\eps^2})
\end{equation}
for some positive constant $C$. Note that here and hereafter $E_x$
is the expectation with respect to the probability law of $X_t$
starting at $X_0=x$.

 If the L\'evy measure part is not absent, the exit
problem  has recently been studied. For symmetric $\alpha$-stable
L\'evy noise, i.e., the L\'evy process whose   L\'evy jump measure
is $\nu(dy) = \frac{dy}{|y|^{1+ \alpha}}$ with $0 < \alpha <2$,
Imkeller and Pavlyukevich \cite{ImkellerP-06} have shown that the
mean exit time is polynomially fast, $O(\frac1{\eps^\a})$, in terms
of $\frac1{\eps}$. Namely, there exist positive constants
$\varepsilon_{0}$, $\gamma$ and $\delta>0$  such that for
$0<\varepsilon \leq \varepsilon_{0}$,
\begin{equation}
E_{x}\sigma(\varepsilon) \sim \frac{\alpha}{\varepsilon^{\alpha}} [
\frac{1}{a^{\alpha}}
+\frac{1}{b^{\alpha}}]^{-1}(1+O(\varepsilon^{\delta})).
\end{equation}
for any $x \in [-b+\gamma, a-\gamma]$.

%Suppose
%\begin{equation}
%Ee^{i \lambda L_{1}}=\exp \{ -d \frac{\lambda^{2}}{2} +\int_{R
%\setminus \{0\} } (e^{ i \lambda y}-1 -i \lambda y I \{|y| <1 \} )
%\frac{dy}{|y|^{1+ \alpha}}\}.
%\end{equation}
%Here $d>0$ and $0<\alpha<2$. That is, $L$ is a sum of a standard
%Brown motion with variance $d>0$ and an independent
%$\alpha$-stable L\'evy motion with $0 < \alpha <2$. The L\'evy measure
%of a stable L\'evy process $L$ is given by
%$v(dy)=\frac{dy}{|y|^{1+\alpha}}$, $y \neq 0$.

Furthermore, for a class of L\'evy noise of exponentially light
jumps, Imkeller, Pavlyukevich and Wetzel \cite{ImkellerP-08} have
shown that the mean exit time is exponentially fast, in terms of
$\frac1{\eps}$, namely,  $O(\exp(\frac1{\eps^\a}))$ with $\a \in (0,
1)$,  or $O(\exp(\frac{|\ln\eps|^{1-\frac1{\a}}}{\eps}) )$ with $\a
>1$.

\vskip 12pt

%\textbf{Theorem } For any $x \in G^{\delta}$, we have
%$lim_{\varepsilon \downarrow 0} E_{x}\varepsilon^{\alpha}
%\sigma(\varepsilon) = \alpha  [ \frac{1}{a^{\alpha}}
%+\frac{1}{b^{\alpha}}]^{-1}$.

In this paper, for a class of non-Gaussian L\'evy noises, we show
that the mean exit time is asymptotically $O(\frac{|\ln
\eps|}{\eps^\a})$. This is   faster than exponential (the well-known
Gaussian Brownian noise case) but slower than the polynomial (the
stable L\'evy noise case).
 So we have   an intermediate  regime,
\begin{equation} \label{compare}
 O(\frac1{\eps^\a} ) <  O(\frac{|\ln \eps|}{\eps^\a})  <
 \exp(\frac{C}{\eps^2}),
\end{equation}
for $0< \eps \ll 1$.

%The method we use is the Large deviation theory and we are to apply
%Godovanchuk's results \cite{God} to study exit problems for
%dynamical systems with L\'evy noise.

\vskip 12pt

In section 2, we   recall the generators for L\'evy processes, and
then prove the main result. In section 3, we consider two examples
of SDEs driven by symmetric L\'evy noises, including  the
$\alpha$-stable symmetric L\'evy noises.

% Freidlin and Wentzell \cite{FW} consider large deviation for locally
%infinitely divisible processes which is a L\'evy process with a jump
%intensity of order $O(\eps^{-1})$ and jump size of order
%$O(\eps)$.More details can be seen in the section 4.

%%%%%%%%%%%%%%%%%%%%%%%%%%%%%%%%%%%%%%%%%%%%%%%%%%%%%%%%%%
\section{Main results }
\label{main}

Let $L_t$ be a L\'evy process with the generating triplet $(\theta, d,
\nu)$.

  It is known that any L\'evy
process is completely determined by the L\'evy-Khintchine formula
(See \cite{Apple, Sato-99, PZ}). This says that for any
one-dimensional L\'evy process $L_{t}$, there exists a $\theta \in R$,
$d>0$ and a measure $\nu$ such that
\begin{equation}
Ee^{i \lambda L_{t}}=\exp \{i \theta\lambda t - d t \frac{\lambda^{2}}{2}
+t \int_{\R \setminus \{0\} } (e^{ i \lambda y}-1 -i \lambda y I
\{|y| <1 \} ) \nu(dy)\},
\end{equation}
where $I(S)$ is the indicator function of the set $S$, i.e., it
takes value $1$ on this set and takes zero value otherwise.

The generator $A$ of the process $L_{t}$ is the same as
infinitesimal generator since L\'evy process has independent and
stationary increments. Hence  $A$ is defined as $A  \phi = \lim_{t
\downarrow 0} \frac{P_{t} \phi -\phi}{t}$ where $P_{t} \phi(x)=
E_{x} \phi(L_{t})$ and $\phi$ is any function belonging to the
domain of the operator $A$. Recall the generator $A$ for  $L_{t}$ is
(See \cite{Apple, PZ})
\begin{equation} \label{A}
A  \phi(x) =  a \phi'(x) +\frac12  d \phi''(x) + \int_{\R
\setminus\{0\}} [\phi(x+ y)-\phi(x) -  \; I\{|y|<1\} \; y \phi'(x) ]
\; \nu(dy).
\end{equation}

Let us   find out the generator of $\eps L_t$.

\begin{lemma} \label{generator}
Let $L_{t}$ be a L\'evy process  with the generating triplet $(a, d,
\nu)$. Then for any $\eps>0$, the generator for $\eps L_t$ is
\begin{equation} \label{Aeps}
A^\eps \phi = \eps a \phi'(x) +\frac12 \eps^2 d \phi''(x) + \int_{\R
\setminus\{0\}} [\phi(x+\eps y)-\phi(x) - \eps  \; I\{|y|<1\} \; y
\phi'(x) ] \; \nu(dy).
\end{equation}
\end{lemma}

\begin{proof}

%%%%%%%%%%%%%%%%%%%%%%%%%%%%%%%%%%%%%%%%%%%%
%Note  that what I used is the observation that with $\eps L_t$, the
%only change in the L\'evy part appears to be in the L\'evy measure which
%is now $\nu (\frac{dy}{\eps})$  So I made a change of variables in
%the above integrals. That is how I get it...... I think I assumed
%somehow the L\'evy measure is "linear" in $y$, i.e., $\nu
%(\frac{dy}{\eps})= \nu (d (\frac{y}{\eps}))$,  which is ok for
%stable L\'evy measure in this section.
%The changes in the drift part
%and Brownian motion part are clear --- We have agreed.
%Note also that this is ``sort of consistent" with eqn (2.1) in \S
%5.2 in \cite{1}, p144. But note that the formulation in  \cite{1}
%may be different from the more modern treatment of L\'evy/Poisson jump
%processes in \cite{Apple, PZ}, for example, in \cite{1}, there is no
%term $I\{  |y|<1\}$. ---I think this question may be partially
%answered from the following proof.

   Notice that
\begin{eqnarray*}
E e^{i \lambda \varepsilon L_{1}} &=&\exp \{i a \eps \lambda -d
\eps^{2} \frac{\lambda^{2}}{2} +\int_{\R \setminus \{0\} } (e^{ i
\lambda \varepsilon y}-1 -i \varepsilon \lambda y I \{|y| <1 \} )
v(dy)\}\\
&=&\exp \{i a \eps \lambda  -i\lambda \eps \int_{\R \setminus
\{0\}}y I \{|y| <1 \} v(dy)  -d \eps^{2} \frac{\lambda^{2}}{2}
+\int_{\R \setminus \{0\} } (e^{ i \lambda y}-1 )
v(d(\frac{y}{\eps}))\}\\
&=&\exp\{i \lambda a \eps -i\lambda \eps\int_{\R \setminus \{0\}}y I
\{|y| <1 \}  v(dy) +i \lambda \int_{\R \setminus \{0\}}y I \{|y|<1\}
v(d(\frac{y}{\eps} )) \\
&\,\,\,\,& -d \eps^{2} \frac{\lambda^{2}}{2}+ \int_{\R \setminus
\{0\} } (e^{ i \lambda y}-1 -i \lambda y I \{|y|<1\} )
v(d(\frac{y}{\eps} ))\}.
\end{eqnarray*}
Hence, $\eps L(t)$ is a L\'evy  process  with the generating triplet
$(\eps a -\eps \int_{\R \setminus \{0\}}y I \{|y| <1 \}
v(dy)+\int_{\R \setminus \{0\}}y I \{|y|<1\} v(d(\frac{y}{\eps}
)),\, d \eps^{2},\, v(d(\frac{y}{\eps})))$. Using the equation
\eqref{A}, it is   seen that the generator of $ \eps L_{t}$ is given
in   \eqref{Aeps}.

This completes the proof of Lemma \ref{generator}.
\end{proof}

\begin{remark}
The generator for the process   $X_{t}^{\eps}$  in \eqref{sde} is
then
\begin{eqnarray} \label{AA}
A^\eps \phi &=& -U^{'}(x) \phi^{'}(x)+ \eps a \phi'(x) +\frac12
\eps^2 d \phi''(x)  \nonumber \\
& +& \int_{\R \setminus\{0\}} [\phi(x+\eps y)-\phi(x) - \eps \; I\{
|y|<1\} \; y \phi'(x) ] \; \nu(dy).
\end{eqnarray}
\end{remark}

\vskip 12pt

%Case 2: If the L\'evy process is light tailed, for example, the L\'evy
%measure of the levy process decreases exponentially fast, I expect
%that exit problem can be studied using action functional under
%appropriate conditions (See theorem 2.1 on page 146 \cite{1}). The
%exit time of the perturbed system $X_{t}^{\varepsilon}$ is expected
%asymptotically exponentially fast.

%Case 3: If L\'evy processes is heavy tailed, exit problem can be
%studied using the large deviation methods using Godovanchuk's
%methods. (see [4] or [8]).

%We first consider the case of stable L\'evy measure and then consider
%more general L\'evy measure.

%Let $I=[-b,a]$, $a,b>0$ and define the first exit time from $I$ as
%\[
%\sigma(\varepsilon) =\inf \{t \geq 0, X_{t}^{\varepsilon} \notin
%[-b,a] \}.
%\]

We make the following assumptions for the SDE \eqref{sde}:

\bigskip

\textbf{(A)} There exists a function $g_{1}(\eps) \downarrow 0$ as
$\eps \downarrow 0$ such that for any $\gamma>0$ $\int_{\R \setminus
[-\gamma,\gamma]}\nu(d(\frac{u}{\eps})) \leq K(\gamma) g_{1}(\eps)$
where $K(\gamma)$ is some function of $\gamma$.

\textbf{(B)} For any $\delta >0$, there exists a positive constant
$K< \infty$ such that $ \int_{\R \setminus
[-K,K]}\nu(d(\frac{u}{\eps})) \leq \delta g_{1}(\eps)$.

\textbf{(C)} There exists a measure $\nu^{*}(du)$ on $\R\setminus
\{0\}$ such that $\frac{1}{g_{1}(\eps)}\nu(d(\frac{u}{\eps}))$
converges weakly to $\nu^{*}(du)$. The limit measure $\nu^{*}$
satisfies the condition that for any Borel set $A \subset \R
\setminus\{0\}$ with measure 0, we have $\nu^{*}(A)=0$.

\textbf{(D)} There exists a function $g_{2}(\eps) \downarrow 0$ as $\eps \downarrow 0$
and a positive constant $K<\infty$ such that
\[
d \eps^{2}+ \int_{\R}\frac{u^{2}}{1+u^{2}}\nu(d
(\frac{u}{\eps}))<Kg_{2}(\eps).
\]
And, there exists some $n>0$ such that $(g_{2}(\eps))^{n}  \leq g_{1}(\eps)$.

\textbf{(E)} There exists a $g_{3}(\eps)\downarrow 0$ as $\eps
\downarrow 0$ and  a positive constant $K<\infty$ such that
\[
\int_{\R}\frac{u}{1+u^{2}}\nu(d(\frac{u}{\eps}))< K g_{3}(\eps).
\]

\medskip

\bigskip

We consider a special class of symmetric L\'evy measures on $\R$ for
$0<\a <2$:
\begin{eqnarray} \label{newlevy}
\nu(du) = f(\ln |u|) \frac{du}{|u|^{1+\a}},
\end{eqnarray}
where $f$ is a nonnegative measurable function on $\R$ such that
this $\nu$ is a L\'evy measure, i.e., it satisfies the above
condition \eqref{levycondition}:
\begin{equation}
  \int_{\R \setminus\{0\} } \frac{u^2}{1+u^2}\; f(\ln |u|) \frac{du}{|u|^{1+\a}}  < \infty.
\end{equation}
Being symmetric, this Levy measure $\nu(du)$ automatically satisfies
the condition (\textbf{E}).

Let $\delta>0$. Define $G^{\delta}=\{x \in [-b,a]: \inf_{t \geq 0}
\min (|Y(0,x,t)-a|, |Y(0,x,t)-(-b)|) \geq \delta\}$.

\begin{theorem} \label{maintheorem}
Consider a class of symmetric L\'evy measures on $\R$  in
  \eqref{newlevy}. Assume that the conditions (\textbf{A}) --
(\textbf{E}) be satisfied. Especially   the condition (\textbf{C})
is satisfied, namely, there exists a positive function
$\tilde{f}(\eps)$ such that
\[ \frac{f(\ln \frac{|u|}{\eps})}{\tilde{f}(\eps)}\frac{du}{|u|^{1+\a}}  \rightharpoonup
\nu^{*}(du).
\]
weakly (as Borel measures on $\R$) to a certain measure $\nu^{*}$.
Then for any $x\in G^\delta$, we have
 \[
 \lim_{\eps \downarrow 0}\tilde{f}(\eps)\;  E_x \;\sigma(\eps)
 = \frac{1}{\nu^{*} (\R \setminus [-b,a])},
 \]
 or  for $\eps \downarrow 0$,
\[
   E_x \;\sigma(\eps)
\sim \frac{1}{\nu^{*} (\R \setminus [-b,a])} \;
\frac1{\tilde{f}(\eps)}.
 \]
\end{theorem}

\begin{proof}
Let us rewrite the generator for the one dimensional Markov
process $X_{t}^{\varepsilon}$.

The generator $A^{\varepsilon}$ depending on parameter
$\varepsilon$ for the one dimensional Markov process
$X_{t}^{\varepsilon}$ in equation \eqref{sde} can be rewritten as
in the following, using Lemma \ref{generator} and \eqref{AA},
\begin{eqnarray*}
A^{\varepsilon}f(x)&=&-U^{'}(x) f^{'}(x)+\eps a f^{'}(x) -\eps
\int_{\R}u I\{|u|<1\} \nu(du) f^{'}(x)  + \frac{d
\varepsilon^{2}}{2}
f^{''}(x) \\
&\,\,&+ \int_{\R} [f(x+ u)-f(x)]
\nu(d(\frac{u}{\eps}))\\
&=&[-U^{'}(x)  +\eps a -\eps  \int_{\R}u I\{|u|<1\} \nu(du)
]f^{'}(x) +
\int_{\R}\frac{u}{1+|u|^{2}}\nu(d(\frac{u}{\eps})) f^{'}(x)\\
&\,\,&  + \frac{d \varepsilon^{2}}{2} f^{''}(x) + \int_{\R}
[f(x+u)-f(x)-\frac{u}{1+|u|^{2}}
f^{'}(x)] \nu(d(\frac{u}{\eps}))\\
&=& U^{\eps}(x) f^{'}(x)+ \frac{d \varepsilon^{2}}{2} f^{''}(x) +
\int_{\R} [f(x+u)-f(x)-\frac{u}{1+|u|^{2}} f^{'}(x)]
\nu(d(\frac{u}{\eps})).
\end{eqnarray*}
Here
\[
U^{\eps}(x)=-U^{'}(x)+\eps a  -\eps  \int_{\R}u I\{|u|<1\}
\nu(du)+\int_{\R}(\frac{u}{1+|u|^{2}})\nu(d(\frac{u}{\eps})).
\]
 The result follows a
similar idea in the proof of the Assertion in  \cite{God}.

\end{proof}

%%%%%%%%%%%%%%%%%%%%%%%%%%%%%%%%%%%%%%%%%%%%%%%%%%%%%%%%%%%%
\section{Examples } \label{example}
%%%%%%%%%%%%%%%%%%%%%%%%%%%%%%%%%%%%%%%%

We look at some applications of the above Theorem \ref{maintheorem}.
\medskip

\textbf{Example 1: $\alpha$-stable symmetric L\'evy Noise}

\medskip
This is a special case of Theorem \ref{maintheorem}  above with
$f(\cdot) \equiv 1$. Consider $X_{t}^{\eps}$ defined in equation
\eqref{sde}, where the L\'evy process $L_{t}$ is  characterized by
\[
Ee^{i \lambda L_{t}}=\exp \{ - t d   \frac{\lambda^{2}}{2} +t
\int_{\R \setminus \{0\} } (e^{ i \lambda u}-1 -i \lambda u I \{|u|
<1 \} ) \frac{1}{|u|^{1+\alpha}}(du)\}.
\]
Here $d \geq 0$ and $0<\alpha <2$ are some constants. The L\'evy
jump measure is $\nu(du)= \frac{1}{|u|^{1+\alpha}}(du)$.  This
 is a so called $\alpha$-stable symmetric L\'evy process,
and it is heavy tailed and has infinite mass due to the strong
intensity of small jumps.

 %Notice that for all $\alpha \in (0,2)$,
%$\int_{R- \{0\}} v(dy)=\int_{R-
%\{0\}}\frac{dy}{|y|^{1+\alpha}}=\infty$ because
%$\int_{0<|y|<1}\frac{dy}{|y|^{1+\alpha}}=\infty$ and $\int_{|y|\geq
%1}\frac{dy}{|y|^{1+\alpha}}<\infty$.  The L\'evy process jumps
%infinitely many times for any given finite time interval   due to
%the strong intensity of small jumps. Furthermore, notice that
%$\int_{R- \{0\}}y^{2} v(dy)=\int_{R-
%\{0\}}y^{2}\frac{dy}{|y|^{1+\alpha}}=\infty$ because
%$\int_{0<|y|<1}y^{2}\frac{dy}{|y|^{1+\alpha}}<\infty$ and
%$\int_{|y|\geq 1}y^{2}\frac{dy}{|y|^{1+\alpha}}=\infty$. This L\'evy
%process is heavy tailed.

Now, we try to verify the conditions in Theorem \ref{maintheorem}.
Notice that
$v(d(\frac{u}{\eps}))=\frac{\varepsilon^{\alpha}}{|u|^{1+\alpha}}du$.
%and for any $\alpha \in (0,2)$,
%\begin{eqnarray*}
%\int_{0}^{\infty} \frac{|u|^{2}}{1+|u|^{2}} v^{\varepsilon}(du)
%&=&\int_{0}^{\infty} \frac{|u|^{2}}{1+|u|^{2}}
%\frac{\varepsilon^{\alpha}}{|u|^{1+\alpha}}
%du \\
%&= & \int_{0}^{1} \frac{|u|^{1-\alpha}}{1+|u|^{2}}
%\varepsilon^{\alpha} du +\int_{1}^{\infty}
%\frac{|u|^{1-\alpha}}{1+|u|^{2}}
%\varepsilon^{\alpha} du\\
%&<& \int_{0}^{1} |u|^{1-\alpha} \varepsilon^{\alpha} du
%+\int_{1}^{\infty} |u|^{-1-\alpha}
%\varepsilon^{\alpha} du\\
%&<& \infty.
%\end{eqnarray*}
It can be    verified that there exist
$g_{1}(\varepsilon)=g_{2}(\varepsilon)=g_{3}(\varepsilon)=\varepsilon^{\alpha}$,
and $\nu^{*}(du)=\frac{1}{|u|^{1+\alpha}}du$ such that the
conditions (\textbf{A})-- (\textbf{D}) are satisfied. Hence, for any
$x\in G^\delta$, we have
 \[
 \lim_{\eps \downarrow 0}\eps^{\alpha}\;  E_x \;\sigma(\eps)
 = \frac{1}{\int_{a}^{\infty}\frac{1}{|u|^{1+\alpha}}du +
 \int_{-\infty}^{-b}\frac{1}{|u|^{1+\alpha}}du}=
 \alpha [\frac{1}{a^{\alpha}}+\frac{1}{b^{\alpha}}]^{-1}.
 \]
 Thus
$$
E_x \;\sigma(\eps) \sim \alpha
[\frac{1}{a^{\alpha}}+\frac{1}{b^{\alpha}}]^{-1} \;
\frac1{\eps^{\alpha}}.
$$
This is the result that Imkeller and Pavlyukevich
\cite{ImkellerP-06} obtained earlier.

\bigskip

\textbf{Example 2: A symmetric L\'evy Noise}
\medskip

Consider $X_{t}^{\eps}$ defined in equation \eqref{sde}, with a
special symmetric   L\'evy process $L_{t}$ that is   characterized
by
\[
Ee^{i \lambda L_{t}}=\exp \{ - d t \frac{\lambda^{2}}{2} +t \int_{\R
\setminus \{0\} } (e^{ i \lambda u}-1 -i \lambda u I \{|u| <1 \} )
\nu(du)\}.
\]
Here $d \geq 0$ and $\nu(du) = f(\ln |u|) \frac{du}{|u|^{1+\a}}$,
$0<\alpha<2$ with $ f(\ln |u|) = \frac1{|\ln |u|| + 1}$. Such a
$\nu(du)$ is a Levy measure satisfying the the condition
\eqref{levycondition}.

We claim that there exist
$g_{1}(\eps)=g_{2}(\eps)=\frac{\eps^{\alpha}}{-\ln \eps}$ and
$\nu^{*}(du)= \frac{1}{u^{1+\alpha}}du$ such that the conditions
(\textbf{A})-(\textbf{D}) are satisfied.

To verify condition (\textbf{A}), it is sufficient to show that for
any $r>0$ there exists some function of $r$, $K(r)$ that
$\int_{r}^{\infty} \nu(d (\frac{u}{\eps})) \leq K(r)
(\frac{\eps^{\alpha}}{-\ln \eps})$ for $\eps$ small enough.  We take
a constant $C=|\min \{\ln r, 0\}|$. Notice that
\begin{eqnarray*}
\int_{r}^{\infty} \nu(d (\frac{u}{\eps})) /
(\frac{\eps^{\alpha}}{-\ln \eps}) &=& \int_{r}^{\infty}
\frac{1}{u^{1+\alpha} (|1+ \frac{ \ln u }{-\ln
\eps}|+\frac{1}{-\ln \eps})}du\\
&\leq& \int_{r}^{\infty} \frac{1}{u^{1+\alpha}[1-\frac{C
}{-\ln \eps}+\frac{1}{-\ln \eps}]} du\\
&=&\int_{r}^{\infty} \frac{1}{u^{1+\alpha}[1+\frac{1-C}{-\ln \eps}]}
du=\frac{1}{\alpha[1+\frac{1-C}{-\ln \eps}]}r^{\alpha}\\
&\leq& \frac{2}{\alpha}r^{\alpha}.
\end{eqnarray*}
for $\eps$ sufficiently small. To verify Condition (\textbf{B}), we
notice that for $K>1 $, $\int_{K}^{\infty} \nu(d (\frac{u}{\eps})) /
(\frac{\eps^{\alpha}}{-\ln \eps}) \leq \int_{K}^{\infty}
\frac{1}{u^{1+\alpha}}du$ which is smaller than any $\delta>0$, if
$K$ is big enough. To verify condition (\textbf{C}), it is
sufficient to show that for any $r>0$,
\[
\lim_{\eps \downarrow 0}
\int_{r}^{\infty}\frac{1}{\frac{\eps^{\alpha}}{\ln \eps}} \nu(d
(\frac{u}{\eps}))= \int_{r}^{\infty} \frac{1}{u^{1+\alpha}}du.
\]
This can be done by a Lebesgue convergence theorem. Finally, let us
verify condition (\textbf{D}). Since $\eps^{2} <
\frac{\eps^{\alpha}}{-\ln \eps}$ for $0<\eps<<1$ and the measure
$\nu$ is symmetric, we only need to show $\int_{0}^{\infty}
\frac{u^{2}}{1+u^{2}} \nu (d \frac{u}{\eps}) /
(\frac{\eps^{\alpha}}{-\ln \eps})$ is bounded by some constant $K$.
Notice that
\begin{eqnarray*}
&\,\,\,\,& \int_{0}^{\infty} \frac{u^{2}}{1+u^{2}} \nu (d
\frac{u}{\eps}) / (\frac{\eps^{\alpha}}{-\ln \eps})\\
&=& \int_{0}^{\infty} (\frac{u^{2}}{1+u^{2}})\frac{1}
{u^{1+\alpha}[|1+\frac{\ln u}{-\ln \eps}| +\frac{1}{-\ln\eps}]}(du) \\
&<& \int_{0}^{ \sqrt{\eps}} \frac{u^{1-\alpha}}{|1+\frac{\ln u}{-\ln
\eps}| +\frac{1}{-\ln\eps}}du+\int_{ \sqrt{\eps}}^{1}
\frac{u^{1-\alpha}}{|1+\frac{\ln u}{-\ln \eps}|
+\frac{1}{-\ln\eps}}du +\int_{1}^{\infty} \frac{1}
{u^{1+\alpha}[|1+\frac{\ln u}{-\ln \eps}| +\frac{1}{-\ln\eps}]}du
\end{eqnarray*}
Here
\[
\int_{0}^{ \sqrt{\eps}} \frac{u^{1-\alpha}}{|1+\frac{\ln u}{-\ln
\eps}| +\frac{1}{-\ln\eps}}du <\int_{0}^{ \sqrt{\eps}}
\frac{u^{1-\alpha}}{\frac{1}{-\ln\eps}}du =\frac{(-\ln\eps) (
\sqrt{\eps})^{2-\alpha}}{2-\alpha} <K_1
\]
for some constant $K_1$ if $\eps$ is sufficiently small. And,
\begin{eqnarray*}
&\,\,\,&\int_{ \sqrt{\eps}}^{1} \frac{u^{1-\alpha}}{|1+\frac{\ln
u}{-\ln \eps}| +\frac{1}{-\ln\eps}}du = \int_{
\sqrt{\eps}}^{1}\frac{u^{1-\alpha}(-\ln \eps)}{|\ln u-\ln \eps|
+1}du\\
&<&\int_{ \sqrt{\eps}}^{1}\frac{u^{1-\alpha}(-\ln
\eps)}{\ln\sqrt{\eps} -\ln \eps +1}du  < \int_{
\sqrt{\eps}}^{1}\frac{u^{1-\alpha}(-\ln \eps)}{
-\frac{1}{2}\ln \eps}du\\
&=& \frac{2}{2-\alpha}[1-(\sqrt{\eps})^{2-\alpha}]<K_2.
\end{eqnarray*}
for some constant $K_2$ if $\eps$ is small enough. Finally,
\[
\int_{1}^{\infty} \frac{1} {u^{1+\alpha}[|1+\frac{\ln u}{-\ln \eps}|
+\frac{1}{-\ln\eps}]}du <\int_{1}^{\infty} \frac{1} {u^{1+\alpha}}du
=\frac{1}{\alpha}.
\]
Therefore, condition (\textbf{D}) is satisfied. By Theorem
\ref{maintheorem}, we conclude that for any $x\in G^\delta$, we have
 \[
 \lim_{\eps \downarrow 0}\frac{\eps^{\alpha}}{-\ln \eps} \;  E_x \;\sigma(\eps)
 = \frac{1}{\int_{a}^{\infty}\frac{1}{|u|^{1+\alpha}}du +
 \int_{-\infty}^{-b}\frac{1}{|u|^{1+\alpha}}du}=
 \alpha [\frac{1}{a^{\alpha}}+\frac{1}{b^{\alpha}}]^{-1}.
 \]
So
$$
E_x \;\sigma(\eps) \sim \alpha
[\frac{1}{a^{\alpha}}+\frac{1}{b^{\alpha}}]^{-1} \;
\frac{|\ln(\eps)|}{\eps^{\alpha}}.
$$

This mean exit time is asymptotically $O(\frac{|\ln
\eps|}{\eps^\a})$. It is   faster than exponential (the well-known
Gaussian Brownian noise case \cite{FW}) but slower than polynomial
(the stable L\'evy noise case \cite{ImkellerP-06}; see also Example
1 above). Namely, for $0< \eps \ll 1$,
\begin{equation} \label{compare2}
 O(\frac1{\eps^\a} ) <  O(\frac{|\ln \eps|}{\eps^\a})  <
 \exp(\frac{C}{\eps^2}).
\end{equation}

\bigskip

{\bf Acknowledgements.} J. Duan would like to thank Professor Ludwig
Arnold for support and encouragement for his  research in random
dynamical systems approach for stochastic systems driven by various
noises.

%%%%%%%%%%%%%%%%%%%%%%%%%%%%%%%%%%%%%%%%%%%%%%%%%%%%%%%%%%%%%%
%%%%%%%%%%%%%%%%%%%%%%%%%%%%%%%%%%%%%%%%%%%%%%%%%%%%%%%%%%%%%%%%%

\end{document}